\documentclass[12pt,twoside]{amsart}
\usepackage{amssymb,times}

\newtheorem{Thm}{Theorem}
\newtheorem{Cor}[Thm]{Corollary}
\newtheorem{Lemma}[Thm]{Lemma}
\newtheorem{Prop}[Thm]{Proposition}
\newtheorem*{Prop*}{Proposition}

\theoremstyle{definition}
\newtheorem{Defn}[Thm]{Definition}
\newtheorem{Ex}{Example}
\newtheorem{Notation}[Thm]{Notation}

\newtheorem{Remark}[Thm]{Remark}

\allowdisplaybreaks[1]

\newcommand{\const}{\text{const}}

\newcommand{\abs}[1]{\left\vert#1\right\vert}
\newcommand{\chf}{\ensuremath{\mathbf{1}}}

\DeclareMathOperator{\E}{E}
\DeclareMathOperator{\Var}{Var}

\newcommand{\set}[1]{\left\{#1\right\}}
\newcommand{\eps}{\varepsilon}

\newcommand{\mf}[1]{\mathbb{#1}}
\newcommand{\mc}[1]{\mathcal{#1}}

\DeclareMathOperator{\Id}{I}
\newcommand{\ip}[2]{\left \langle #1, #2 \right \rangle}

\DeclareMathOperator{\Res}{\mc{R}}
\DeclareMathOperator{\Gen}{\mc{D}}
\DeclareMathOperator{\esssup}{\mathrm{esssup}}

\title{Free martingale polynomials}
\author[M.~Anshelevich]{Michael Anshelevich}
\thanks{This work was supported in part by an NSF postdoctoral fellowship, and by an MSRI postdoctoral fellowship}
\address{Department of Mathematics, University of California, Berkeley, CA 94720}
\email{manshel@math.berkeley.edu}
\date{\today}

\begin{document}

\begin{abstract}
In this paper we investigate the properties of the free Sheffer systems, which are certain families of martingale polynomials with respect to the free L\'{e}vy processes. First, we classify such families that consist of orthogonal polynomials; these are the free analogs of the Meixner systems. Next, we show that the fluctuations around free convolution semigroups have as principal directions the polynomials whose \emph{derivatives} are martingale polynomials. Finally, we indicate how Rota's finite operator calculus can be modified for the free context.
\end{abstract}

\maketitle

\section{Introduction}

Hermite polynomials $H_n(x,t) = (-1)^n t^n e^{x^2/2t} \partial_x^n e^{-x^2/2t}$ are related to the Gaussian convolution semigroup $\mu_t(d x) = \frac{1}{\sqrt{2 \pi t}} e^{-x^2/2t} d x$ in two ways. First, for every fixed $t$, the polynomials $\set{H_n(x,t)}_{n=0}^\infty$ are the monic orthogonal polynomials with respect to $\mu_t$. Second, they are martingale polynomials for the corresponding L\'{e}vy process, namely the Brownian motion. This means that if $\set{B(t)}$ is a Brownian motion, for each $n$ the process $H_n(B(t),t)$ is a martingale with respect to the standard filtration of $\set{B(t)}$. This easily follows from the fact that the exponential generating function of these polynomials $\sum_{n=0}^\infty \frac{1}{n!} H_n(x,t) z^n$ is precisely $e^{xz - t z^2/2}$, the exponential martingale for the Brownian motion. This result goes back at least to \cite{McKean}.

Polynomials whose exponential generating function has a general form of this type $f(z)^t e^{u(z) x}$ are called Sheffer systems. These systems, especially a particular sub-class of them called the Appell systems, have been investigated in depth (see, for example, \cite{Lai} and their references, as well as the references of our Section \ref{Sec:OperatorCalculus}). If in addition the polynomials are orthogonal, they are called Meixner systems. There is a complete classification of these systems due to Meixner, described in detail in \cite[Chapter 4]{SchOrthogonal}. See also \cite[Chapter 5]{FeiSch} for the description of the same objects from a somewhat different point of view. Up to various re-scalings, Meixner systems form a one-parameter family consisting of the Meixner / Laguerre / Meixner-Pollaczek polynomials, with the Hermite and Charlier polynomials obtained as limiting cases.

The first part of this paper is concerned with the investigation of the corresponding objects in \emph{free probability}. This is a non-commutative probability theory, in which the usual independence is replaced by a different notion of \emph{free independence}, and the usual convolution is replaced by the notion of (additive) \emph{free convolution}. Over the last twenty years, this theory has exhibited depth which may someday rival that of the classical probability theory. It also exhibits surprising analogy with the usual probability theory; the structure underlying this analogy still remains largely a mystery. The results of this paper provide further examples of this analogy. Namely, we define the free Sheffer systems to be the systems of polynomials which are martingales for processes with freely independent increments; see the precise definitions in the next section. Then the free Meixner systems are the free Sheffer systems consisting of orthogonal polynomials. It was known that the Chebyshev polynomials of the 2nd kind are martingale polynomials for the free Brownian motion \cite{BiaMultiplicativeBM}, and it follows from the results of \cite{AnsFSM} that the corresponding statement holds for the free Charlier polynomials and the free Poisson process. We show here that up to re-scaling, the free Meixner systems also form a one-parameter family, with the aforementioned free analogs of the Hermite and the Charlier systems  arising as limits. The free Meixner polynomials are much simpler than the classical ones: their recursion relations have almost constant coefficients. Nevertheless, the analogy with the classical case is exact. This is especially surprising since the free Meixner systems do not correspond to the classical ones in the canonical bijection between the classical and the free infinitely divisible distributions; see Section~\ref{Subsection:Classical}.

The second part of this paper was motivated by the article \cite{CDFluctuations}. There, using free stochastic calculus, Cabanal-Duvillard re-proved the result of Johannson \cite{Johansson} that the principal directions for the fluctuations around the semicircular limit for large Gaussian random matrices are given by the Chebyshev polynomials of the first kind; he also re-proved the corresponding result for the Wishart matrices, and extended both results to pairs of random matrices. The method of proof involves precisely the kind of martingale and orthogonality properties we are considering here. Instead of the random matrix context, in this paper we consider a semigroup of operators of convolution with a family of freely infinitely divisible distributions. These operators are non-linear, and we consider their differentials. We show that the principal directions for these differentials are given by polynomials whose derivatives are martingale polynomials.

In the third part of the paper, we begin the investigation of the free Sheffer systems using the finite operator calculus machinery of Rota. The original finite operator calculus describes precisely the classical Sheffer systems. Multiple generalizations of that calculus have been considered; in particular, free binomial sequences fit into one of such generalizations. Interestingly, however, the free Sheffer sequences can also be considered in the context of a different kind of finite operator calculus, which apparently has not been investigated before.

\noindent\textbf{Acknowledgments.} I would like to thank Franz Lehner, Jim Pitman, Thierry Cabanal-Duvillard, and Dan Voiculescu for useful conversations.

\section{Preliminaries}

\subsection{Formal power series}

Let $\mf{C}[x]$ be the algebra of complex polynomials in an indeterminate $x$. We will frequently abuse notation by denoting, for example, by $x^n$ the function $x \mapsto x^n$. Denote by $C^1(\mf{R}_+)$ the space of complex-valued differentiable functions, and by $C^1(\mf{R}_+)[x]$ the space of polynomials in $x$ with coefficients that are differentiable functions of $t \in \mf{R}_+$. We will consider formal power series $H(t,z) = \sum_{k=0}^\infty a_k(t) z^{-k}$ in $z^{-1}$ and formal Laurent series $H(t,z) = \sum_{k = -n}^\infty a_k(t) z^k$ in $z$. More generally, we will consider formal power series $H(x,t,z) = \sum_{k=0}^\infty P_n(x,t) z^n$ in $z$, where $P_n(x,t)$ is a polynomial in $x$ of degree $n$ with $t$-dependent coefficients. For a formal power series $H(z)$, $\frac{1}{H(z)}$ will always denote its inverse under multiplication; this is well-defined iff $H(0) \neq 0$. $H^{-1}(z)$ will always denote its inverse under composition; this is well-defined iff $H(0) = 0$, $H'(0) \neq 0$. Denote by $\mc{P}_{0,1}$ the space of all formal power series $u$ in $z$ with coefficients in $\mf{C}$ such that $u(0) = 0$, $u'(0) = 1$.

For a linear functional $\nu$ on $\mf{C}[x]$, denote $\ip{\nu}{p}$ the value of $\nu$ on $p \in \mf{C}[x]$. Denote $m_n(\nu) = \ip{\nu}{x^n}$, the $n$-th moment of $\nu$. Clearly the functional $\nu$ can be identified with its moment sequence $\set{m_n(\nu)}_{n=0}^\infty$. Denote by $\mc{M}$ the space of all linear functionals on $\mf{C}[x]$, by $\mc{M}_1$ the subset $\set{\nu \in \mc{M} | \ip{\nu}{1} = 1}$ of unital functionals, and by $\mc{M}_0 = \set{\nu \in \mc{M} | \ip{\nu}{1} = 0}$ the orthogonal complement to the constants.

\subsection{Difference quotient}

Define the canonical derivation $\partial: \mf{C}[x]  \rightarrow \mf{C}[x] \otimes \mf{C}[x]$ by the requirement that $\partial(1) = 0$, $\partial(x) = 1 \otimes 1$. If we identify $\mf{C}[x] \otimes \mf{C}[x]$ with $\mf{C}[x,y]$, then $\partial f(x,y) = \frac{f(x)-f(y)}{x-y}$, the difference quotient. Moreover, define the maps $\partial^k: \mf{C}[x] \rightarrow \mf{C}[x]^{\otimes(k+1)}$ by $\partial^k = k (1 \otimes \cdots \otimes 1 \otimes \partial) \partial^{k-1}$. More explicitly, on monomials their action is 
\begin{equation*}
\partial^k x^n = k! \sum_{\substack{i(0), i(1), \ldots, i(k) \geq 0 \\ i(0) + i(1) + \ldots + i(k) = n-k}} x^{i(0)} \otimes x^{i(1)} \otimes \cdots \otimes x^{i(k)}.
\end{equation*}
Note that if $m$ is the multiplication map $\mf{C}[x] \otimes \mf{C}[x] \rightarrow \mf{C}[x]$, then $m \circ \partial = \partial_x$.

For a formal power series $H(x,t,z)$, by $\partial_t, \partial_x, \partial$, or $\nu$ applied to it we mean the formal power series in $z$ obtained by the term-wise application of these operations. In particular, for $z \in \mf{C} \backslash \mf{R}$, $x \in \mf{R}$, denote by $\Res_z$ the resolvent function $x \mapsto \frac{1}{z-x} = \sum_{n=0}^\infty x^n z^{-(n+1)}$. It has the property that $\partial \Res_z = \Res_z \otimes \Res_z$, and more generally $\partial^k \Res_z = k! \Res_z^{\otimes (k+1)}$.

\subsection{Free convolution machinery}

For the background in free probability, the reader should consult the main references \cite{VDN,VoiSF}; whenever we don't give specific references one of these can be used. See also \cite{SpeNCReview} for a survey of the combinatorial approach to free probability.

Let $\nu \in \mc{M}$. Define the Cauchy transform of $\nu$ to be the formal power series in $1/z$, $G_\nu(z) = \ip{\nu}{\Res_z} = \sum_{n=0}^\infty m_n(\nu) z^{-(n+1)}$. If $\nu \in \mc{M}_1$, $m_0(\nu) = 1$. Therefore for such $\nu$, the Cauchy transform series has an inverse under composition of the form $K_\nu(z) = \frac{1}{z} + \sum_{n=1}^\infty r_n(\nu) z^n$, where $\set{r_n(\nu)}$ are the \emph{free cumulants} of $\nu$. Finally, denote $R_\nu(z) = K_\nu(z) - \frac{1}{z}$ the \emph{$R$-transform} of $\nu$, which is a power series in $z$. The main property of the $R$-transform is that $R_{\mu \boxplus \nu} = R_\mu+ R_\nu$, where $\boxplus$ is the operation of additive free convolution, which in this paper will be called simply free convolution. This is a certain commutative, associative, non-linear binary operation on probability measures, which can be extended to an operation on $\mc{M}_1$. The above property of the $R$-transform can be taken as the definition of $\boxplus$. See the references for its relation to free independence, and also to the lattice of noncrossing partitions.

Let $\mu$ be a probability measure all of whose moments are finite. Then to it naturally corresponds an element of $\mc{M}_1$, although this correspondence is neither injective nor surjective. In particular, the above notions apply to it. In fact, in this case $G_\mu(z) = \int_{\mf{R}} \frac{1}{z - x} d \mu(x)$ is an analytic function on $\mf{C} \backslash \mf{R}$, and $K_\mu, R_\mu$ are analytic functions on a Stolz angle in $\mf{C}_+$. In this case $\mu$ can be recovered from its Cauchy transform $G_\mu$ by taking a weak limit:
\begin{equation}
\label{WeakLimit}
\mu(dx) = - \frac{1}{\pi} \lim_{\eps \rightarrow 0^+} \Im G_\mu(x + i \eps) d x.
\end{equation}

From now on we assume that $\mu$ is an freely infinitely divisible distribution. This means that there exists a free convolution semigroup $\set{\mu_t}_{t \in [0, \infty)}$ of probability measures, characterized by the properties that $\mu_0 = \delta_0$, $\mu_t \boxplus \mu_s = \mu_{t+s}$, $\mu_1 = \mu$, $R_{\mu_t} = t R_\mu$. Throughout the paper $G_{\mu_t}, K_{\mu_t}, R_{\mu_t}$ will be denoted by, respectively, $G_t, K_t, R_t$. Denote $F_{s,t} = K_s \circ G_t$. It satisfies $F_{s,t}(0) = 0$, $F_{s,t}'(0) = 1$. Let $\nu \in \mc{M}_1$, $\nu_t = \nu \boxplus \mu_t$, and $G(t,z) = G_{\nu_t}(z)$. Then these formal power series satisfy a quasi-linear differential equation
\begin{equation}
\label{Quasilinear}
\partial_t G(t,z) + R_\mu(G(t,z)) \partial_z G(t,z) = 0.
\end{equation}

Throughout the paper, it should be clear from the context (whether the objects treated are measures or functionals) whether we are considering analytic functions or only formal power series.

\subsection{Noncommutative stochastic processes}

Let $(\mc{A}, \E)$ be a noncommutative probability space. That is, $\mc{A}$ is a finite von Neumann algebra, $\E$ is a faithful normal tracial state on $\mc{A}$, and $\tilde{\mc{A}}$ is the algebra of unbounded operators affiliated to $\mc{A}$. Let $\set{X(t)}_{t \in [0, \infty)}$ be a free L\'{e}vy process on $\mc{A}$ with distribution $\set{\mu_t}$. That is, for all $t$, $X(t)$ is a self-adjoint operator in $\tilde{\mc{A}}$, $X(0) = 0$, the distribution of $X(t)$ with respect to $\E$ is $\mu_t$, and for any $0 < t_1 < t_2 < \ldots < t_n$, the family $\set{X_{t_1}, X_{t_2} - X_{t_1}, \ldots, X_{t_n} - X_{t_{n-1}}}$ is a freely independent family. Let $\set{\mc{A}_t}$ be the natural filtration of $\mc{A}$ induced by the process $\set{X(t)}$, and let $\E_t: \mc{A} \rightarrow \mc{A}_t$ be the trace-preserving conditional expectations. They are characterized by the property that for $X \in \mc{A}, Y \in \mc{A}_t, \E[X Y] = \E[\E_t[X] Y]$.

It was proven in \cite{BiaProcesses} that the process $\set{X(t)}$ is a Markov process with respect to the natural filtration $\set{\mc{A}_t}$. More specifically, for any bounded Borel function $f$,
\[
\E_s[f(X(t))] = (\mc{K}_{s,t}(f))(X(s)).
\]
Here $\set{\mc{K}_{s,t}}_{0 \leq s \leq t}$ is a family of Feller integral operators, characterized by the property that
\[
\mc{K}_{s,t}(\Res_z) = \Res_{F_{s,t}(z)}
\]
for any $z \in \mf{C} \backslash \mf{R}$. Since under our assumptions, $F_{s,t}$ has a formal power series expansion, it follows that the operators $\mc{K}_{s,t}$ can be extended to operators on $\mf{C}[x]$.

\subsection{Martingale polynomials}

Let a function $p(x,t)$, $\mf{R} \times \mf{R}_+ \rightarrow \mf{C}$ be, for each $t$, bounded and measurable in $x$. We will call it a martingale function for $\set{\mu_t}$ if for all $s < t$, $\mc{K}_{s,t}(p(\cdot, t)) = p(\cdot, s)$. Then using the above characterization, $p$ is a martingale function if and only if the process $t \mapsto p(X(t), t)$ is a martingale, i.e. an $\set{\mc{A}_t}$-measurable process such that for $s<t$, $\E_s[p(X(t), t)] = p(X(s), s)$. In particular, for any $z \in \mf{C} \backslash \mf{R}$, the function $\frac{1}{z(K_t(z) - x)}$ is a martingale function. These are the analogs of the exponential martingales for the usual L\'{e}vy processes. More generally, let $\Omega$ be a domain in $\mf{C} \backslash \mf{R}$, and $u, v$ be functions on it such that $0 \not \in u(\Omega)$, $\forall t \in \mf{R}_+, K_t(v(\Omega)) \subset \mf{C} \backslash \mf{R}$. Then for $z \in \Omega$, the process
\begin{equation}
\label{Martingale}
t \mapsto \frac{1}{u(z)(K_t(v(z)) - X(t))}
\end{equation}
is also a martingale. If $p(x,t)$ is a polynomial in $x$ such that for all $s < t$, $\mc{K}_{s,t}(p(\cdot, t)) = p(\cdot, s)$, we call it a martingale polynomial for the semigroup $\set{\mu_t}$.

Since all of the moments  of $\mu$, and hence of all $\mu_t$, are finite, the function $K_t$ has a power series expansion 
\[
K_t(z) = \frac{1}{z} + t \sum_{n=1}^\infty r_n z^{n-1},
\]
where $\set{r_n}$ are the free cumulants of $\mu$. Suppose the functions $u$ and $v$ have formal power series expansions so that $u, v \in \mc{P}_{0,1}$. Then we can define the polynomials $\set{Q_n(x,t)}_{n=0}^\infty$ by their generating function
\begin{equation*}
\frac{1}{u(z)(K_t(v(z)) - x)} = H(x,t,z) = \sum_{n=0}^\infty Q_n(x,t) z^n.
\end{equation*}
Note that $Q_n(x,t)$ has degree $n$ as a polynomial in $x$, and its highest coefficient is equal to $1$, independently of $t$. Moreover, it is also a polynomial of degree $n$ in $t$. Since there is an open set $\Omega$ such that for $z \in \Omega$, $H(x,t,z)$ is well-defined and so its power series expansion converges, it follows from \eqref{Martingale} that for all $n$, the process $t \mapsto Q_n(X(t), t)$ is also a martingale. We will call any such family of martingale polynomials a generalized free Sheffer system for $\set{\mu_t}$. If $u = v$, we will call it a \emph{free Sheffer system}. Finally, if $u(z) = v(z) = z$, we will call it a standard Sheffer system for $\set{\mu_t}$; the term \emph{free Appell system} would also be appropriate.

\begin{Lemma}
\label{Lem:StandardSheffer}
Any martingale polynomial for $\set{\mu_t}$ is a linear combination of the elements of the standard Sheffer system for it.
\end{Lemma}

\begin{proof}
Any martingale polynomial has a constant highest term coefficient (since $\E_s[X(t)^n] = \E_s[(X(s) + (X(t) - X(s)))^n] = X(s)^n + $ lower order terms). So any martingale polynomial of degree $n$ is a linear combination of the element of the standard Sheffer system of degree $n$ and a martingale polynomial of degree at most $(n-1)$. The result follows by induction.
\end{proof}

\subsection{Orthogonal polynomials}

Since $\mu$ has moments of all orders, we can define $\set{P_n(x,t)}_{n=0}^\infty$ to be the family of monic polynomials orthogonal with respect to $\set{\mu_t}$ (by which we mean that for each $t \in \mf{R}_+$, $\set{P_n(\cdot, t)}$ are orthogonal with respect to $\mu_t$). They will satisfy $P_0(x,t) = 1$ and a $3$-term recursion relation
\[
P_{n+1}(x,t) = (x - \alpha_{n+1}(t)) P_n(x,t) - \beta_n(t) P_{n-1}(x,t)
\]
for $n \geq 0$, with the convention that $P_{-1} = 0$, and all $\beta_n(t) \geq 0$. We will denote $\gamma_n = \ip{\mu_t}{P^2_n(\cdot, t)} = \prod_{j=1}^n \beta_j$.

In general the polynomials $\set{P_n(x,t)}$ will not be martingale polynomials.

\begin{Defn}
A family of polynomials $\set{P_n(x,t)}$ orthogonal with respect to a free convolution semigroup $\set{\mu_t}$ that is also a generalized free Sheffer system for that semigroup is a \emph{free Meixner system}.
\end{Defn}

\subsection{Semicircular distributions}

The semicircular distribution with mean $\alpha$ and variance $\beta$ is
\[
\sigma_{\alpha, \beta}(d x) = \frac{1}{2 \pi \beta} \sqrt{4 \beta - (x - \alpha)^2} \chf_{[\alpha - 2 \sqrt{\beta}, \alpha + 2 \sqrt{\beta}]}(x) d x.
\]
Denote $\sigma_t = \sigma_{0, t}$. Then $\set{\sigma_t}$ is a free convolution semigroup. See the main references and also the beginning of Section \ref{Sec:Fluctuations} for its importance. Also, the arcsine distribution with mean $\alpha$ and variance $2 \beta$ is
\[
\sigma_{\alpha, \beta}'(d x) = \frac{1}{\pi} \frac{1}{\sqrt{4 \beta - (x - \alpha)^2}} \chf_{[\alpha - 2 \sqrt{\beta}, \alpha + 2 \sqrt{\beta}]}(x) d x.
\]
The monic Chebyshev polynomials of the second kind $\set{U_n(x,t)}$ are the orthogonal polynomials with respect to $\set{\sigma_t}$; they are defined by $U_n(x,t) = t^{n/2} U_n(x/\sqrt{t})$, $U_n(2 \cos \theta) = \frac{\sin (n+1) \theta}{\sin \theta}$. The Chebyshev polynomials of the first kind are the orthogonal polynomials with respect to $\set{\sigma_{0, t}'}$; they are defined by $T_n(x,t) = t^{n/2} T_n(x/\sqrt{t})$, $T_n(2 \cos \theta) = \cos n \theta$. Both families satisfy the recursion relations $P_{n+1}(x) = x P_n(x) - P_{n-1}(x)$, with initial conditions $U_0(x) = 1, U_1(x) = x$, $T_0(x) = 1, T_1(x) = \frac{1}{2} x$. They are related by $\partial_x T_n(x,t) = \frac{n}{2} U_{n-1}(x,t)$.

\section{Free Meixner systems}

\begin{Lemma}
\label{Lemma:Sheffer}
Let $\nu$ be a probability measure all of whose moments are finite. Let $\set{P_n}_{n=0}^\infty$ be the family of monic polynomials orthogonal with respect to $\nu$, satisfying
\[
P_{n+1}(x) = (x - \alpha_{n+1}) P_n(x) - \beta_n P_{n-1}(x)
\]
for $n \geq 0$, with all $\beta_n \geq 0$. 
\begin{enumerate}
\item 
\label{condition}
The generating function $H(x,z) = \sum_{n=0}^\infty P_n(x) z^n$ has the form 
\[
H(x,z) = \frac{1}{u(z) (f(z) - x)}
\]
for $u, u f$ having formal power series expansions with $u \in \mc{P}_{0,1}$, $(u f)(0) = 1$ if and only if for some $\alpha, \alpha', \beta, \beta'$ with $\beta, \beta - \beta' \geq 0$, $\alpha_n = \alpha - \delta_{n1} \alpha'$, $\beta_n = \beta - \delta_{n1} \beta'$.
\item 
In this case $f = K_\nu \circ u$. 
\item 
Let 
\[
Q_n(x) = U_n(x - \alpha, \beta).
\]
Then under the conditions of \eqref{condition},
\begin{align}
\label{PQ}
P_0(x) &= Q_0(x), \notag \\
P_1(x) &= Q_1(x) + \alpha' Q_0(x), \\
P_n(x) &= Q_n(x) + \alpha' Q_{n-1}(x) + \beta' Q_{n-2}(x) \notag
\end{align}
for $n > 1$.
\end{enumerate}
\end{Lemma}

Note that the polynomials $\set{Q_n}$ above are orthogonal with respect to the semicircular distribution $\sigma_{\alpha, \beta}$ of mean $\alpha$ and variance $\beta$.

\begin{proof}
First suppose that $H$ is of the above form. The polynomials $\set{P_n}$ form a basis of $\mf{C}[x]$, so we may define the lowering operator $A$ on $\mf{C}[x]$ by $A P_n = P_{n-1}$ for $n \geq 0$, and extend linearly. Then
\[
A(H)(x,z) = \sum_{n=0}^\infty (A P_n)(x) z^n = \sum_{n=0}^\infty P_{n-1}(x) z^n = z H(x,z).
\]
On the other hand, 
\[
x H(x,z) = f(z) H(x,z) - \frac{1}{u(z)}.
\]
Therefore $A(x H)(x,z) = z f(z) H(x,z) = z x H(x,z) + \frac{z}{u(z)}$. The second term has a formal power series expansion, $\frac{z}{u(z)} = \sum_{n=0}^\infty c_n z^n$. Thus finally, $A(x P_n)(x) = x P_{n-1}(x) + c_n$. Now apply the operator $A$ to the recursion relation. We obtain
\[
P_n(x) = (x - \alpha_{n+1}) P_{n-1}(x) - \beta_n P_{n-2}(x) + c_n
\]
for $n \geq 1$. Subtracting from it the recursion relation for $n-1$, we obtain
\[
(\alpha_n - \alpha_{n+1}) P_{n-1}(x) + (\beta_{n-1} - \beta_n) P_{n-2}(x) + c_n = 0
\]
for $n \geq 1$. The polynomials $\set{P_n}$ are linearly independent for different $n$. We conclude that $\alpha_1 - \alpha_2 + c_1 = 0$, $\alpha_n - \alpha_{n+1} = 0$ for $n > 1$, $\beta_1 - \beta_2 + c_2 = 0$, $\beta_{n-1} - \beta_n = 0$ for $n > 2$. Therefore the recursion relations in fact have the form 
\begin{align*}
P_0(x) &= 1, \\
P_1(x) &= x - (\alpha - \alpha'), \\ 
P_2(x) &= (x - \alpha) P_1(x) - (\beta - \beta') P_0(x), \\
P_{n+1}(x) &= (x - \alpha) P_n(x) - \beta P_{n-1}(x)
\end{align*}
for $n \geq 2$, for some $\alpha, \alpha', \beta, \beta'$ with $\beta \geq 0$, $\beta - \beta' \geq 0$.

Conversely, for polynomials with such recursion relations, 
\[
H(x,z) 
= \frac{1 + \alpha' z + \beta' z^2}{1 + \alpha z + \beta z^2 - x z} 
= \frac{1}{u(z) (f(z) - x)},
\]
with
\[
u(z) = \frac{z}{1 + \alpha' z + \beta' z^2}
\]
and
\[
f(z) = \frac{1}{z} + \alpha + \beta z.
\]
Since $\ip{\nu}{P_n(\cdot) P_0(\cdot)} = \delta_{n 0}$,
\[
1 
= \ip{\nu}{H(\cdot,z)} 
= \frac{1}{u(z)} \ip{\nu}{\Res_{f(z)}}
= \frac{1}{u(z)} G_\nu(f(z)),
\]
and so $u(z) = G_\nu(f(z))$, $f = K_\nu \circ u$.

The expression for the polynomials $\set{P_n}$ in terms of the polynomials $\set{Q_n}$ follows from the fact that the latter satisfy the recursion relations 
\[
Q_{n+1}(x) = (x - \alpha) Q_n(x) - \beta Q_{n-1}(x)
\]
for $n \geq 0$. 
\end{proof}

\begin{Remark}
Orthogonal polynomials with constant recursion coefficients have been described in \cite{CTConstant}. The argument with the lowering operator above is similar to the original one of Meixner as described in \cite{SchOrthogonal}; see also Section \ref{Sec:OperatorCalculus}. Finally, for the free Poisson case (see below) the description of the orthogonal polynomials in terms of the shifted Chebyshev polynomials has appeared in \cite{HagThoRMExact}.
\end{Remark}

In the following theorem, the cases are labeled by the names of the distributions and the orthogonal polynomials for the corresponding classical Meixner systems.

\begin{Thm}
\label{Thm:Meixner}
Up to affine transformations of $x$ and re-scaling $t$ by a positive factor, the following are all the non-trivial free Meixner systems. We list the recursion relations for the polynomials, the corresponding free convolution semigroup, its $R$-transform, the function $u$ such that
\[
\sum_{n=0}^\infty P_n(x,t) z^n 
= H(x,t,z) 
= \frac{1}{u(z)(K_t(u(z)) - x)} 
= \frac{1}{1 + t u(z) R_\mu(u(z)) - u(z) x},
\]
and the expression in terms of shifted Chebyshev polynomials of the 2nd kind.
\begin{description}
\item[Semicircular / Chebyshev] This case corresponds to the classical Gaussian / Hermite case.
\[
P_{n+1}(x,t) = x P_n(x,t) - t P_{n-1}(x,t)
\]
for $n \geq 0$. 
\begin{gather*}
P_n(x,t) = U_n(x, t), \\
\mu_t(d x) = \sigma_t(d x), \\
u(z) = z, \qquad R_\mu(z) = z.
\end{gather*}
\item[Poisson / Charlier] $P_1(x,t) = x - t$,
\[
P_{n+1}(x,t) = (x - (t + 1)) P_n(x,t) - t P_{n-1}(t)
\]
for $n \geq 1$.
\[
P_n(x,t) = Q_n(x,t) + Q_{n-1}(x,t)
\]
for $n \geq 1$, where $Q_n(x,t) = U_n(x - (t+1), t)$.
\begin{gather*}
\mu_t(d x) =
\frac{t}{x} \sigma_{1+t, t}(d x) + \max(1-t, 0) \delta_0, \\
u(z) = \frac{z}{1+z}, \qquad R_\mu(z) = \frac{1}{1-z}.
\end{gather*}

In the remaining three cases, for a parameter $a \geq 0$,
\begin{align*}
P_1(x,t) &= x - a t, \\
P_2(x,t) &= (x - a (t+2)) P_1(x,t) - t P_0(x,t), \\
P_{n+1}(x,t) &= (x - a (t + 2)) P_n(x,t) - (t + 1) P_{n-1}(t)
\end{align*}
for $n \geq 2$; thus $\alpha(t) = a (t+2)$, $\beta(t) = t+1$.
\begin{align*}
P_1(x,t) &= Q_1(x,t) + 2 a Q_0(x,t), \\
P_n(x,t) &= Q_n(x,t) + 2 a Q_{n-1}(x,t) + Q_{n-2}(x,t)
\end{align*}
for $n \geq 2$, where $Q_{n}(x,t) = U(x - \alpha(t), \beta(t))$. Also,
\[
u(z) = \frac{z}{1 + 2a z + z^2}
\]
and
\[
R_\mu(z) 
= \frac{z^{-1} - \sqrt{(2a - z^{-1})^2 - 4}}{2}.
\]
The free convolution semigroups are as follows.
\item[Continuous binomial / Meixner-Pollaczek] $0 \leq a < 1$.
\[
\mu_t(d x)
= \frac{t \beta(t)}{x^2 + t^2 (1 - a^2)} \sigma_{\alpha(t), \beta(t)}(d x).
\]
\item[Gamma / Laguerre] $a = 1$.
\[
\mu_t(d x) = \frac{t (1+t)}{x^2} \sigma_{2+t, 1+t}(d x).
\]
\item[Negative binomial / Meixner] $a > 1$.
\[
\mu_t(d x)
= \frac{t \beta(t)}{x^2 - t^2 (a^2 - 1)} \sigma_{\alpha(t), \beta(t)}(d x) + \max \Bigl(1 - t \frac{a - \sqrt{a^2-1}}{2 \sqrt{a^2-1}}, 0 \Bigr) \delta_{t \sqrt{a^2 - 1}}(x).
\]
\end{description}
\end{Thm}

Note that we make no claim that these families exhaust the situations when the orthogonal polynomials with respect to a free convolution semigroup $\set{\mu_t}$ are also martingale functions for it: we restrict the analysis to the generating functions of a specific form. However, see Lemma \ref{Lem:Generator}.

\begin{proof}
Let $\set{P_n(x,t)}$ be a free Meixner system, and let $\set{\mu_t}$ be the corresponding free convolution semigroup. By Lemma \ref{Lemma:Sheffer}, $\alpha_n(t) = \alpha(t) - \delta_{n1} \alpha'(t)$, $\beta_n(t) = \beta(t) - \delta_{n1} \beta'(t)$, and 
\[
H(x,t,z) = \frac{1}{u(t,z)(f(t,z) - x)},
\]
for 
\[
u(t,z) = \frac{z}{1 + \alpha'(t) z + \beta'(t) z^2}
\]
and 
\[
f(t,z) = \frac{1}{z} + \alpha(t) + \beta(t) z.
\]
For a free Meixner system, $u$ does not depend on $t$, and 
\[
f(t,z) = K_t(u(z)) = \frac{1}{u(z)} + t R(u(z)).
\]
In particular, $\set{P_n(x,t)}$ is in fact a free Sheffer system, rather than a generalized one. We also conclude that $\alpha'(t) = \alpha'$, $\beta'(t) = \beta'$, $\alpha(t) = a_1 t + a_2$, $\beta(t) = b_1 t + b_2$, with $b_1, b_2 \geq 0$. Moreover, since the measures $\set{\mu_t}$ form a free convolution semigroup, the expectation and the variance of $\mu_t$ are proportional to $t$. For the measure $\mu_t$, its expectation is equal to $\alpha(t) - \alpha'$ and its variance is equal to $\beta(t) - \beta'$. Therefore $\alpha' = a_2, \beta' = b_2$. The case $b_1 = 0$ is a degenerate case of zero variance, so assume $b_1 > 0$, and in fact re-normalize $t$ so that $b_1 = \Var(\mu_1) = 1$. We conclude that
\[
u(z) = \frac{z}{1 + a_2 z + b_2 z^2},
\]
and
\begin{equation}
\label{R}
R_\mu \left( \frac{z}{1 + a_2 z + b_2 z^2} \right) = a_1 + b_1 z.
\end{equation}
Let 
\[
w = \frac{z}{1 + a_2 z + b_2 z^2},
\]
i.e.
\begin{equation}
\label{w}
b_2 z^2 + (a_2 - w^{-1})z + 1 = 0. 
\end{equation}
First suppose $b_2 = 0$. Then 
\[
z = \frac{1}{w^{-1} - a_2} = \frac{w}{1 - a_2 w},
\]
and so
\[
R_\mu(w) = a_1 + \frac{w}{1 - a_2 w}.
\]
By adding a constant to $x$ we may assume that $a_1 = a_2$. For $a_2 = 0$ we obtain the semicircular distribution. For $a_2 \neq 0$, we may re-scale $x$ so that $a_1 = a_2 = 1$. We obtain the free Poisson distribution. See the more complicated cases below for the method.

From now on, assume $b_2 \neq 0$. By re-scaling $x$ and $t$ we may assume that $b_2 = b_1 = 1$. By adding a constant to $x$ we may assume that $a_2 = 2 a_1$, and denote $a_1$ by $a$. Possibly by replacing $P_n(x,t)$ by $(-1)^n P_n(-x, t)$ we may assume that $a \geq 0$. Thus the recursion relation takes the form 
\[
P_{n+1}(x,t) = (x - a (t + 2)) P_n(x,t) - (t + 1) P_{n-1}(x,t)
\]
for $n \geq 2$, with $\alpha(t) = a (t+2)$, $\beta(t) = t+1$.

From equations~\eqref{R} and~\eqref{w},
\begin{gather*}
R_\mu(w) 
= a + \frac{-(2a - w^{-1}) - \sqrt{(2a - w^{-1})^2 - 4}}{2} 
= \frac{w^{-1} - \sqrt{(2a - w^{-1})^2 - 4}}{2}, \\
K_t(w) 
= \frac{(2 + t) w^{-1} - t \sqrt{(2a - w^{-1})^2 - 4}}{2},
\end{gather*}
so
\begin{equation}
\label{Cauchy}
G_t(z) = \frac{(2 + t) z - t^2 a - t \sqrt{(z - \alpha(t))^2 - 4 \beta(t)}}{2(z^2 - t^2 (a^2 - 1))}.
\end{equation}
Using equation \eqref{WeakLimit}, the formulas for $\mu_t$ follow from equation \eqref{Cauchy}. Note that in the negative binomial case, there is at most one atom, and there are no atoms in the gamma case $a=1$.

The expressions for the orthogonal polynomials in terms of the Chebyshev polynomials of the 2nd kind follow immediately from Lemma~\ref{Lemma:Sheffer}.
\end{proof}

\begin{Remark}
\label{Remark:Meixner}
More generally, we could consider the situation when the family of polynomials $\set{P_n}$ is a generalized free Sheffer system for the free convolution semigroup $\set{\mu_t}$ but is orthogonal with respect to some other family of measures $\set{\nu_t}$ which do not form a free convolution semigroup. In this case we can conclude that $\alpha(t) = a_1 t + a_2$, $\beta(t) = b_1 t + b_2$,
\begin{gather*}
v(z) = \frac{z}{1 + a_2 z + b_2 z^2}, \\
u(z) = \frac{z}{1 + \alpha' z + \beta' z^2}, \\
R_\mu \left( \frac{z}{1 + a_2 z + b_2 z^2} \right) = a_1 + b_1 z,
\end{gather*}
and 
\[
K_{\nu_t} (u(z)) = K_t(v(z)) = f(t,z) = \frac{1}{z} + \alpha(t) + \beta(t) z.
\]
Therefore the measures $\set{\mu_t}$ form one of the families in Theorem~\ref{Thm:Meixner}, and the measures $\set{\nu_t}$ are measures of the same type, but possibly with different parameters and scaling. In particular, whenever such a family $\set{\nu_t}$ exists, the orthogonal polynomials with respect to the family $\set{\mu_t}$ themselves form a free Meixner system. We do not calculate the measures $\set{\nu_t}$ explicitly, except for one case: $\alpha' = \beta' = 0$, $u(z) = z$. In this case $\nu_t$ is the semicircular distribution with mean $\alpha(t)$ and variance $\beta(t)$, and the orthogonal polynomials are 
\[
Q_{n}(x,t) = U_n(x - \alpha(t), \beta(t)).
\]
Since the relations \eqref{PQ} (with $a_2$ in place of $\alpha'$, $b_2$ in place of $\beta'$) are invertible, we have already observed that these are martingale polynomials for $\set{\mu_t}$.
\end{Remark}

\subsection{Relation to classical orthogonal polynomials}
\label{Subsection:Classical}

There is a bijective correspondence between the classical and the free infinitely divisible measures investigated in detail in \cite{BerPatDomains}; see \cite{AnsQCum} for a simple description of it in the case of measures all of whose moments are finite. This bijection naturally transforms limit theorems for independent random variables into limit theorems for freely independent random variables. In particular, it maps the normal distribution to the semicircular distribution, and the Poisson distribution to the free Poisson distribution (hence the name). We now show that, surprisingly, except for these two cases, the correspondence provided by Theorem~\ref{Thm:Meixner}, that maps a classical Meixner system to the free Meixner system with the same parameter $a$, is different from this bijection.

For measures with finite variance, the Bercovici-Pata bijection $\Lambda$ takes the following form. Let $\nu$ be an infinitely divisible measure with mean $\lambda$ and canonical measure $\tau$. That is, denoting by $\mc{F}_\nu$ the Fourier transform of $\nu$, $(\log \mc{F}_\nu)'(0) = \lambda$ and $(\log \mc{F}_\nu)''(\theta) = \mc{F}_\tau(\theta)$. Then $\Lambda(\nu)$ is the freely infinitely divisible measure with mean $\lambda$ and free canonical measure $\tau$. That is, $R_{\Lambda(\nu)}(z^{-1}) = \lambda + G_\tau(z)$.

For a free Meixner system with parameter $a$, $R_\mu(0) = a$. Thus
\[
R_\mu(w^{-1}) = a + \frac{w - 2a - \sqrt{(2a - w)^2 - 4}}{2}.
\]
Therefore the free canonical measure of the free convolution semigroup corresponding to the parameter $a$ is the semicircular distribution with mean $2a$ and variance $1$.

On the other hand, for $a=1$ the classical polynomials are the Laguerre polynomials. They are orthogonal with respect to the standard gamma distribution, which has the Fourier transform $(1 - i \theta)^{-1}$. But $(\log (1 - i \theta)^{-1})'' = (1 - i \theta)^{-2}$. Therefore its canonical measure is the gamma distribution with parameter $2$ (in both cases the scaling parameter of the gamma distribution is taken to be $1$). Thus, $\Lambda$ does not map the classical Meixner system with parameter $1$ to a free Meixner system.

\begin{Remark}[$q$-interpolation]
\label{Rem:QDeform}
A possible explanation for the correspondence in Theorem \ref{Thm:Meixner} is provided by the following interpolating family of polynomials. For $n \geq 0$ and $q \in [-1, \infty)$, denote by $[n]_q = \sum_{i=0}^{n-1} q^{i}$ the $q$-integer, with the convention that $[0]_q = 0$. Define the $q$-Hermite polynomials by the recursion relation
\[
P_{n+1}(x,t) = x P_n(x,t) - t [n]_q P_{n-1}(x,t),
\]
the $q$-Charlier polynomials by
\[
P_{n+1}(x,t) = (x - (t + [n]_q)) P_n(x,t) - t [n]_q P_{n-1}(x,t),
\]
and the $q$-Meixner family with parameter $a$ by
\[
P_{n+1}(x,t) = (x - a (t + 2 [n]_q)) P_n(x,t) - [n]_q(t + [n-1]_q) P_{n-1}(x,t),
\]
all for $n \geq 0$. Then we obtain the classical families for $q=1$, and the free families for $q=0$. See \cite{SaiKraw,AnsQCum}, Proposition \ref{Prop:General}, Lemma~\ref{Lemma:Counter}, and the discussion following it for further results about these families. 
\end{Remark}

\begin{Remark}[IID-Sheffer systems]
If the parameter $t$ is discrete rather than continuous, there are many more families of martingale polynomials. A standard example are the Krawtchouk polynomials $P_n(x,N)$. Given a parameter $p \in [0, 1]$, these are orthogonal with respect to the corresponding binomial distribution, which is a convolution of $N$ copies of the Bernoulli distribution $(1-p) \delta_0 + p \delta_1$. In free probability, the corresponding measure is the distribution of the sum of $N$ freely independent projections. This distribution is easy to find, and in \cite{FreeKraw} the corresponding orthogonal polynomials, the free Krawtchouk polynomials, have been calculated; see also \cite{SaiKraw}. For the sum of $N$ freely independent projections of trace $p$, the corresponding orthogonal polynomials satisfy the recursion relations
\begin{align*}
P_1(x,N) &= x - N p, \\
P_2(x,N) &= (x - (N p + (1 - 2p)) P_1(x,N) - N p (1-p) P_0(x,N), \\
P_{n+1}(x,N) &= (x - (N p + (1 - 2p)) P_n(x,N) - (N - 1) p (1-p) P_{n-1}(x,N)
\end{align*}
for $n \geq 2$. Thus, $\alpha(N) = a_1 N + a_2$, $\beta(N) = b_1 N + b_2$, with $a_1 = p$, $a_2 = 1 - 2p$, $b_1 = p (1-p)$, $b_2 = -p (1-p)$. Note that $b_2 < 0$, so these recursion relations are not of the standard free Meixner form. The corresponding generating function is
\begin{align*}
H(x,N,z) 
&= \frac{1 + (1-2p) z - p(1-p) z^2}{(1 + (1-2p) z - p(1-p) z^2 ) + N (p z + p(1-p) z^2) - z x} \\
&= \frac{1}{1 + N u(z) R(u(z)) - u(z) x},
\end{align*}
with
\[
u(z) = \frac{z}{1 + (1-2p) z - p(1-p) z^2}
\]
and 
\[
R(z) 
= \frac{1 - z^{-1} - \sqrt{((1-2p) - z^{-1})^2 + 4p(1-p)}}{2}
= R_{(1-p) \delta_0 + p \delta_1}(z).
\]
This implies that the free Krawtchouk polynomials are martingale polynomials for the free binomial process. That is, let $\set{p_i}_{i=1}^\infty$ is a family of freely independent projections of trace $p$, and let $X(N) = \sum_{i=1}^N p_i$, for $N \geq 1$. Then for the corresponding free Krawtchouk polynomials $\set{P_n(x,N)}$ and $N_0 < N$, the conditional expectation of $P_n(X(N), N)$ onto the von Neumann algebra generated by $\set{X(k)}_{k=1}^{N_0}$ is $P_n(X(N_0), N_0)$.
\end{Remark}

\section{Fluctuations}
\label{Sec:Fluctuations}

\begin{Notation}
A family $\set{P_n}_{n=0}^\infty$ of polynomials such that the degree of $P_n$ is $n$ form a basis in $\mf{C}[x]$. Denote by $\set{P_n^\ast}$ the dual basis of $\mc{M}$, determined by $\ip{P_n^\ast}{P_k} = \delta_{nk}$. If the polynomials $\set{P_n}$ are orthogonal with respect to a probability measure $\nu$ that is uniquely determined by its moments, they are an orthogonal basis for $L^2(\mf{R}, \nu)$, and so $\set{P_n^\ast}$ is a basis for the dual space of measures $L^2(\mf{R}, \nu)'$. In this case we can identify explicitly that $P_n^\ast(d x) = \frac{1}{\gamma_n} P_n(x) \nu(d x)$, where $\gamma_n = \ip{\nu}{P_n^2}$.
\end{Notation}

Let $\tau$ be a probability measure with mean $0$ and variance $1$. Denote by $S_c$ the scaling operator, $S_c(\tau)(\Omega) = \tau(c^{-1} \Omega)$. Then the free central limit theorem states that $\tau^{\boxplus n} \circ S_{1/\sqrt{n}} \rightarrow \sigma_1$ weakly. If we call $C$ the operator $\tau \mapsto (\tau \boxplus \tau) \circ S_{1/\sqrt{2}}$, then $\sigma_1$ is a fixed point of $C$, and the theorem says that it is an attracting fixed point, $C^n \tau \rightarrow \sigma_1$. In \cite{AnsCLT}, we investigated the fluctuations around this limit, and showed that the derivative $D_{\sigma_1} C$ of $C$ at $\sigma_1$ has eigenfunctions $T_n^\ast$ with eigenvalues $2^{1 - n/2}$, where $\set{T_n}$ are the Chebyshev polynomials of the first kind. Here by the derivative $D_\tau C(\nu)$ we will mean the G\^{a}teaux derivative
\[
\lim_{\eps \rightarrow 0} \frac{1}{\eps}(C(\tau + \eps \nu) - C(\tau))
\]
when the limit exists in the appropriate topology. In a remark in \cite{CDFluctuations}, Cabanal-Duvillard re-interpreted this result as follows:
\[
\sigma_t \boxplus (\sigma_s + \eps T_n^\ast(s)) 
= \sigma_{t+s} + \eps T_n^\ast(t+s) + o(\eps).
\]
Given this,
\begin{align*}
((\sigma_1 + \eps T_n^\ast) \boxplus (\sigma_1 + \eps T_n^\ast)) \circ S_{1/\sqrt{2}}
&= \sigma_2 \circ S_{1/\sqrt{2}} + 2 \eps T_n^\ast(2) \circ S_{1/\sqrt{2}} + o(\eps) \\
&= \sigma_1 + 2 \eps 2^{-n/2} T_n^\ast \circ S_{\sqrt{2}} \circ S_{1/\sqrt{2}} \\
&= \sigma_1 + \eps 2^{1 - n/2} T_n^\ast,
\end{align*}
so the previous result follows.

In this section we extend this analysis to all free convolution semigroups (Corollary~\ref{Cor:Fluctuation}) and, in a more precise sense, to all free Meixner systems (Corollary~\ref{Cor:Meixner}).

\begin{Notation}
Let $\mu \in \mc{M}_1$. Denote by $C_\mu$ the operator of free convolution with $\mu$ on $\mc{M}_1$, $C_\mu: \nu \mapsto \mu \boxplus \nu$. For $\set{\mu_t}$ a free convolution semigroup, denote $C_{\mu_t}$ simply by $C_t$. Note that unlike in the classical case, $C_\mu$ is a non-linear operator.

$\mc{M}$ is a topological vector space with the weak$^\ast$ topology induced on it as a dual space of $\mf{C}[x]$. Note that if a sequence $\set{\tau_n}$ of elements of $\mc{M}$ correspond to measures, and if the limit of this sequence corresponds to a \emph{unique} measure, then the corresponding sequence of measures converges weakly.

The tangent space to $\mc{M}_1$ at any point is naturally identified with $\mc{M}_0$. So for $\mu, \tau \in \mc{M}_1$, $\nu \in \mc{M}_0$, we can define the G\^{a}teaux derivative $D_\tau C_\mu(\nu) = \lim_{\eps \rightarrow 0} \frac{1}{\eps}(C_\mu(\tau + \eps \nu) - C_\mu(\tau))$ when the limit exists in the above  topology on $\mc{M}_0$. This limit always exists; we delay the proof of this fact until Lemma~\ref{Lemma:Composition}.
\end{Notation}

\begin{Lemma}
\label{Lem:Generator}
Let $\set{\mu_t}$ be a free convolution semigroup, and $\nu \in \mc{M}_1$. Denote $\nu_t = \nu \boxplus \mu_t$. Let $\set{\Gen_{t, \nu}}_{t \in \mf{R}_+}$ be the family of operators on $C^1(\mf{R}_+)[x]$ given by
\begin{equation}
\label{Generator}
\Gen_{t, \nu} f = \sum_{n=1}^\infty \frac{1}{(n-1)!} r_n \ip{\Id \otimes \nu_t^{\otimes (n-1)}}{(\partial_x \otimes \Id^{\otimes (n-1)}) \partial^{n-1} f} + \partial_t f.
\end{equation}
\begin{enumerate}
\item
\label{all}
For any $f \in C^1(\mf{R}_+)[x]$, $\partial_t \ip{\nu_t}{f} = \ip{\nu_t}{\Gen_{t, \nu} f}$.
\item
\label{mart}
For a martingale polynomial $f \in C^1(\mf{R}_+)[x]$ and $\nu = \delta_0$, $\Gen_{t, \delta_0} f = 0$.
\end{enumerate}
\end{Lemma}

\begin{proof}
It suffices to prove the first part for a monomial $f(x,t) = a(t) x^n$. Denoting $m_n(t) = \ip{\nu_t}{x^n}$,
\[
\partial_t \ip{\nu_t}{f} = a'(t) m_n(t) + a(t) m_n'(t) = \ip{\nu_t}{\partial_t f} + a(t) \partial_t \ip{\nu_t}{x^n}.
\]
Therefore it suffices to prove the statement for all the monomials $x^n$. Let
\[
\Res_z(x) = \sum_{n=0}^\infty x^n z^{-(n+1)}
\]
be their formal generating function. Then $G(t,z) = \ip{\nu_t}{\Res_z} = \sum_{n=0}^\infty m_n(t) z^{-(n+1)}$. On the other hand,
\begin{equation}
\label{DRes}
\begin{split}
\Gen_{t, \nu} \Res_z
&= \sum_{n=1}^\infty \frac{1}{(n-1)!} r_n \ip{\Id \otimes \nu_t^{\otimes (n-1)}}{(\partial_x \otimes \Id^{\otimes (n-1)}) \partial^{n-1} \Res_z} + \partial_t \Res_z \\
&= \sum_{n=1}^\infty r_n \ip{\Id \otimes \nu_t^{\otimes (n-1)}}{(\partial_x \otimes \Id^{\otimes (n-1)}) \Res_z^{\otimes n}} \\
&= \sum_{n=1}^\infty r_n G(t,z)^{n-1} \Res_z^2
= - R_\mu(G(t,z)) \partial_z \Res_z.
\end{split}
\end{equation}
So for $\Res_z$, the desired equation takes the form
\[
\partial_t G(t,z) = - R_\mu(G(t,z)) \partial_z G(t,z),
\]
which is equation~\eqref{Quasilinear}.

Now we consider the second part. By Lemma~\ref{Lem:StandardSheffer}, it suffices to show the second property for the polynomials in the standard Sheffer system, or indeed for their generating function $H(x,t,z) = \frac{1}{z (K_t(z) - x)}$. But
\[
\Gen_{t, \delta_0} H(x,t,z) = \sum_{n=1}^\infty r_n \frac{1}{z(K_t(z) - x)^2} G_t(K_t(z))^{n-1} - \frac{R_\mu(z)}{z(K_t(z) - x)^2}
= 0. 
\]
\end{proof}

\begin{Remark}
Suppose that in the preceding lemma, $\nu$ is in fact a probability measure. Then there is a free L\'{e}vy process $\set{X_t}$ and an operator $Y_0$ freely independent from it, so that the distribution of $Y_t = Y_0 + X_t$ is $\nu_t$. Then 
\begin{multline}
\label{ExpectationGenerator}
(\mc{D}_{t, \nu} f)(Y(t), t) = \lim_{h \rightarrow 0} \frac{1}{h} \left( \E_t[f(Y(t + h), t + h)] - f(Y(t), t) \right) \\
= \partial_h \bigl|_{h = 0} \E_t[f(Y(t + h), t + h)]
\end{multline}
and
\[
(\mc{D}_{t, \nu} f)(x,t) = \lim_{h \rightarrow 0} \frac{1}{h} \left( (\mc{K}_{t,t + h}f)(x,t + h) - f(x,t) \right) = \partial_h \bigl|_{h = 0} (\mc{K}_{t,t + h}f)(x,t + h).
\]
That is, $\set{\mc{D}_{t, \nu}}$ are the generators of the family of operators $\set{\mc{K}_{s,t}}$. In this case, the results of the preceding lemma follow immediately from equation~\eqref{ExpectationGenerator}. The expression~\eqref{Generator} for $\mc{D}_{t, \nu}$ follows essentially from equation~\eqref{DRes}: 
\begin{align*}
\mc{D}_{t, \nu} \Res_z 
&= \partial_h \bigl|_{h = 0} \Res_{F_{t , t + h}(z)}
= \partial_h \bigl|_{h = 0} \frac{1}{K_t(G_{t + h}(z)) - x}
= \partial_h \bigl|_{h = 0} \frac{1}{z - h R_\mu(G_{t + h}(z)) - x} \\
&= \frac{R_\mu(G_t(z))}{(z - x)^2}
= \sum_{n=1}^\infty \frac{1}{(n-1)!} r_n \ip{\Id \otimes \nu_t^{\otimes(n-1)}}{(\partial_x \otimes \Id^{\otimes(n-1)}) \partial^{n-1} \Res_z}.
\end{align*}
For compactly supported $\set{\mu_t}$, this conclusion also follows from the functional It\^{o} formula for the free L\'{e}vy processes obtained in \cite{AnsIto}.
\end{Remark}

\begin{Prop}
Let $p(\cdot, t)$ be a martingale polynomial for $\set{\mu_t}$. Let $g(\cdot, t)$ be a function such that $\partial_x g(x,t) = p(x,t)$. Let $\nu \in \mc{M}_0$. Then 
\[
\ip{D_{\mu_s} C_t (\nu)}{g(\cdot, t+s)} = \ip{\nu}{g(\cdot, s)}.
\]
\end{Prop}

\begin{proof}
We will show that
\[
\ip{\mu_t \boxplus (\mu_s + \eps \nu)}{g(\cdot, t+s)} = \ip{\mu_{t+s}}{g(\cdot, t+s)} + \eps \ip{\nu}{g(\cdot, s)} + o(\eps). 
\]
Denote $\nu_t = \mu_t \boxplus (\mu_s + \eps \nu)$. Apply part \eqref{mart} of Lemma~\ref{Lem:Generator} to $p(\cdot, t) = \partial_x g(\cdot,t)$.
\begin{align*}
0 
&= \sum_{n=1}^\infty \frac{1}{(n-1)!} r_n \ip{\Id \otimes \mu_t^{\otimes (n-1)}}{(\partial_x \otimes \Id^{\otimes (n-1)}) \partial^{n-1}  \partial_x g(\cdot,t)} + \partial_t \partial_x g(\cdot,t) \\
&= \sum_{n=1}^\infty \frac{1}{(n-1)!} r_n \ip{\Id \otimes \mu_t^{\otimes (n-1)}}{ \Bigl( (\partial_x^2 \otimes \Id^{\otimes (n-1)}) + (n-1) (\partial_x \otimes \Id^{\otimes (n-2)} \otimes \partial_x) \Bigr) \partial^{n-1} g(\cdot,t)} \\*
&\quad + \partial_x \partial_t g(\cdot,t) \\
&= \partial_x \biggl( \sum_{n=1}^\infty \frac{1}{(n-1)!} r_n \ip{\Id \otimes \mu_t^{\otimes (n-1)}}{ \Bigl( (\partial_x \otimes \Id^{\otimes (n-1)}) + (n-1) (\Id^{\otimes (n-1)} \otimes \partial_x) \Bigr) \partial^{n-1} g(\cdot,t)} \\*
&\quad + \partial_t g(\cdot,t) \biggr) \\
&= \partial_x \biggl( \sum_{n=1}^\infty \frac{1}{(n-1)!} r_n \ip{\Id \otimes \mu_t^{\otimes (n-1)} + (n-1) \mu_t^{\otimes (n-1)} \otimes \Id}{(\partial_x \otimes \Id^{\otimes (n-1)}) \partial^{n-1} g(\cdot,t)} \\*
&\quad + \partial_t g(\cdot,t) \biggr).
\end{align*}
That is,
\[
\sum_{n=1}^\infty \frac{1}{(n-1)!} r_n \ip{\Id \otimes \mu_t^{\otimes (n-1)} + (n-1) \mu_t^{\otimes (n-1)} \otimes \Id}{(\partial_x \otimes \Id^{\otimes (n-1)}) \partial^{n-1} g(\cdot,t)} + \partial_t g(\cdot,t) = \const
\]
as a function of $x$. Now apply part \eqref{all} of Lemma \ref{Lem:Generator}, and expand $(\mu_t \boxplus (\mu_s + \eps \nu))^{\otimes n}$ in powers of $\eps$: 
\[
(\mu_t \boxplus (\mu_s + \eps \nu))^{\otimes n} 
= \mu_{t+s}^{\otimes n} + \eps \sum_{i=1}^n \mu_{t+s}^{\otimes (i-1)} \otimes (\nu_t - \mu_{t+s}) \otimes \mu_{t+s}^{\otimes (n-i)} + o(\eps).
\]
We obtain
\begin{align*}
&\partial_t \ip{\nu_t}{g(t+s)} \\
&\quad = \sum_{n=1}^\infty \frac{1}{(n-1)!} r_n \ip{\nu_t^{\otimes n}}{(\partial_x \otimes \Id^{\otimes (n-1)}) \partial^{n-1} g(\cdot,t+s)} + \ip{\nu_t}{\partial_t g(\cdot,t+s)} \\
&\quad = \sum_{n=1}^\infty \frac{1}{(n-1)!} r_n \ip{\mu_{t+s}^{\otimes n}}{(\partial_x \otimes \Id^{\otimes (n-1)}) \partial^{n-1} g(\cdot,t+s)} + \ip{\mu_{t+s}}{\partial_t g(\cdot,t+s)} \\*
&\quad \quad +  \eps \sum_{n=1}^\infty \frac{1}{(n-1)!} r_n \Bigl\langle ((\nu_t - \mu_{t+s}) \otimes \mu_{t+s}^{\otimes (n-1)}) + (n-1) (\mu_{t+s}^{\otimes (n-1)} \otimes (\nu_t - \mu_{t+s})), \\*
&\quad \quad \quad (\partial_x \otimes \Id^{\otimes (n-1)}) \partial^{n-1} g(\cdot,t+s) \Bigr\rangle + \eps \ip{\nu_t - \mu_{t+s}}{\partial_t g(\cdot,t+s)} + o(\eps) \\
&\quad = \partial_t \ip{\mu_{t+s}}{g(\cdot,t+s)} \\*
&\quad \quad +  \eps (\nu_t - \mu_{t+s}) \biggl( \sum_{n=1}^\infty \frac{1}{(n-1)!} r_n \Bigl\langle (\Id \otimes \mu_{t+s}^{\otimes (n-1)}) + (n-1) (\mu_{t+s}^{\otimes (n-1)} \otimes \Id), \\*
&\quad \quad \quad (\partial_x \otimes \Id^{\otimes (n-1)}) \partial^{n-1} g(\cdot,t+s) \Bigr\rangle + \partial_t g(\cdot,t+s) \biggr) + o(\eps) \\
&\quad = \partial_t \ip{\mu_{t+s}}{g(\cdot,t+s)} + o(\eps),
\end{align*}
since $\ip{\nu_t - \mu_{t+s}}{1} = 0$. Also,
\[
\ip{\nu_0}{g(\cdot,s)} 
= \ip{\mu_s + \eps \nu}{g(\cdot,s)}
= \ip{\mu_s}{g(\cdot,s)} + \eps \ip{\nu}{g(\cdot,s)}.
\]
Therefore $\ip{\nu_t}{g(t+s)} = \ip{\mu_{t+s}}{g(\cdot,t+s)} + \eps \ip{\nu}{g(s)} + o(\eps)$.
\end{proof}

\begin{Cor}
\label{Cor:Fluctuation}
Let $\set{V_n(x,t)}_{n=1}^\infty$ be a family of \emph{fluctuation polynomials} for $\set{\mu_t}$, that is, any family such that for $n \geq 1$, $\partial_x V_n(x,t)$ is a martingale polynomial for $\set{\mu_t}$ of degree $n-1$. Then for each $t \in \mf{R}_+$, $\set{1} \cup \set{V_n(\cdot, t)}_{n=1}^\infty$ is a basis for $\mf{C}[x]$. Denote by $\set{V_n^\ast(t)}_{n=1}^\infty$ the dual basis of $\mc{M}_0$. Then for $n \geq 1$, $D_{\mu_s} C_t(V_n^\ast(s)) = V_n^\ast(t+s)$.
\end{Cor}

\begin{proof}
The function $\partial_x V_n(x,t)$ is a martingale polynomial.
\begin{align*}
\ip{\mu_t \boxplus (\mu_s + \eps V_n^\ast(s))}{V_k(\cdot,t+s)}
&= \ip{\mu_{t+s}}{V_k(\cdot,t+s)} + \eps \delta_{nk} + o(\eps) \\
&= \ip{\mu_{t+s} + \eps V_n^\ast(t+s)}{V_k(\cdot,t+s)} + o(\eps).
\end{align*}
Since they also take the same value on the scalar $1$, we conclude that $\mu_t \boxplus (\mu_s + \eps V_n^\ast(s)) = \mu_{t+s} + \eps V_n^\ast(t+s) + o(\eps)$.
\end{proof}

\begin{Remark}[Algebraic infinite divisibility]
\label{Rem:AllInfinitelyDivis} 
For \emph{any} $\mu \in \mc{M}_1$, we can define a free convolution semigroup $\set{\mu_t} \subset \mc{M}_1$ by $R_{\mu_t} = t R_\mu$; the distinguishing characteristic of infinitely divisible probability measures is that for them all of $\set{\mu_t}$ are in fact positive measures. But the analysis of this section, in particular the preceding corollary, applies equally well to such purely algebraic semigroups.
\end{Remark}

\begin{Cor}
\label{Cor:Meixner}
Let $\set{\mu_t}$ be a free convolution semigroup associated to a free Meixner system, and let $\alpha(t), \beta(t)$ be the parameters of the corresponding recursion relations. Then the polynomials $V_n(x,t) = T_n(x - \alpha(t), \beta(t))$ are a family of fluctuation polynomials for $\set{\mu_t}$. In particular, $D_{\mu_s} C_t: L^2_0(\mf{R}, \sigma_{\alpha(s), \beta(s)}')' \rightarrow L^2_0(\mf{R}, \sigma_{\alpha(t+s), \beta(t+s)}')'$. Here $L^2_0(\mf{R}, \tau)' = \set{\nu \in L^2(\mf{R}, \tau)' | \ip{\tau}{1} = 0}$.
\end{Cor}

\begin{Remark}[Classical case]
Throughout this remark only, let $\set{\mu_t}$ be a convolution semigroup with respect to the usual convolution $\ast$, and let $C_t$ be the operator of the usual convolution with $\mu_t$. Then $C_t$ itself is a linear operator. Let $\set{P_n(x,t)}$ be a collection of martingale polynomials for the corresponding classical L\'{e}vy process. Then it is easy to see that for the dual basis to these polynomials themselves, $C_t P_n^\ast(s) = P_n^\ast(s+t)$. Note that the standard Sheffer system of polynomials has a generating function $e^{x z - \log \mc{F}(z,t)}$, whose derivative is $z e^{x z - \log \mc{F}(z,t)}$. Thus in the classical case, a derivative of a martingale polynomial is again a martingale polynomial. In particular, for a classical Meixner system, $C_t: L^2_0(\mf{R}, \mu_s)' \rightarrow L^2_0(\mf{R}, \mu_{t+s})'$.

In fact, in the classical case, such a statement holds for any convolution operator $C_\mu$. We want to show that for Borel probability measures $\mu, \sigma$, if $\nu \in L^2(\mf{R}, \sigma)'$, then $\nu \ast \mu \in L^2(\mf{R}, \sigma \ast \mu)'$. By definition, $\nu \in L^2(\mf{R}, \sigma)'$ iff $\forall f \in L^2(\mf{R}, \sigma)$, $\ip{\nu}{f} < \infty$. Therefore, it suffices to show that for $f \in L^2(\mf{R}, \sigma \ast \mu)$, $\ip{\nu \ast \mu}{f} < \infty$. But
\[
\ip{\mu \ast \nu}{f}
= \ip{\nu}{f \ast \check{\mu}},
\]
where $d \check{\mu}(x) = d \mu(-x)$. So it suffices to show that $f \ast \check{\mu} \in L^2(\mf{R}, \sigma)$, in other words, that the operator of convolution with $\check{\mu}$ maps $L^2(\mf{R}, \sigma \ast \mu)$ into $L^2(\mf{R}, \sigma)$.

For $f \in L^1(\mf{R}, \sigma \ast \mu)$, 
\begin{align*}
\int \abs{(f \ast \check{\mu})(x)} d \sigma(x)
&= \int \abs{\int f(x) d \mu(x - y)} d \nu(y) \\
&\leq \iint \abs{f(x)} d \mu(x - y) d \nu(y)
= \int \abs{f(x)} d (\mu \ast \sigma)(x) 
< \infty.
\end{align*}
For $f \in L^\infty(\mf{R}, \sigma \ast \mu)$, 
\[
\esssup \abs{f \ast \check{\mu}} \leq \esssup \abs{f} < \infty.
\]
The result follows by Riesz interpolation.
\end{Remark}

\subsection{Cauchy transforms}

Let $\set{P_n(x,t)}$ be the standard Sheffer system for $\set{\mu_t}$. Then
\[
\mc{K}_{s,t} (P_n(\cdot, t)) = P_n(\cdot, s).
\]
Therefore the adjoint operator $\mc{K}^\ast_{s,t}$ on $\mc{M}$ is defined and determined by 
\[
\ip{\mc{K}^\ast_{s,t} (P_n^\ast(s))}{P_k(\cdot, t)}
= \ip{P_n^\ast(s)}{\mc{K}_{s,t} (P_k(\cdot, t))}
= \ip{P_n^\ast(s)}{P_k(\cdot, s)}
= \delta_{nk}.
\]
That is, $\mc{K}^\ast_{s,t} (P_n^\ast(s)) = P_n^\ast(t)$.

\begin{Lemma}
For $\nu \in \mc{M}$, $G_{\mc{K}_{s,t}^\ast \nu}(z) = G_\nu(F_{s,t}(z))$. In particular, if $\set{P_n(x,t)}$ is a family of martingale polynomials for $\set{\mu_t}$ and $\set{P_n^\ast(t)}$ is the dual basis of $\mc{M}$, then $G_{P_n^\ast(t)}(z) = G_{P_n^\ast(s)}(F_{s,t}(z))$.
\end{Lemma}

\begin{proof}
We only need to prove the first statement.
\[
G_{\mc{K}_{s,t}^\ast \nu}(z) 
= \ip{\mc{K}^\ast_{s,t} \nu}{\Res_z} 
= \ip{\nu}{\Res_{F_{s,t}(z)}}
= G_{\nu}(F_{s,t}(z)).
\]
\end{proof}

\begin{Ex}
For a free Meixner system $\set{P_n(x,t)}$, we can find $G_{P_n^\ast(t)}$ explicitly. If $H(x,t,z)$ is the generating function of the polynomials $\set{P_n(x,t)}$,
\[
\frac{1}{G_t(z) (z - x)} = H(x, t, u^{-1}(G_t(z))) = \sum_{n=0}^\infty (u^{-1}(G_t(z)))^n P_n(x,t).
\]
Thus
\[
(u^{-1}(G_t(z)))^n = \ip{P^\ast_n(t)}{H(\cdot, t, u^{-1}(G_t(z)))} = \frac{1}{G_t(z)} G_{P^\ast_n(t)}(z),
\]
so
\[
G_{P^\ast_n(t)}(z) = G_t(z) (u^{-1}(G_t(z)))^n.
\]
For the modified Chebyshev polynomials of Remark~\ref{Remark:Meixner}, we obtain by the same method
\begin{equation}
\label{GChebyshev2}
G_{Q^\ast_n(t)}(z) = u(G_t(z)) G_t(z)^n.
\end{equation}
\end{Ex}

\begin{Lemma}
\label{Lemma:Composition}
Let $\nu \in \mc{M}_0$.
\begin{enumerate}
\item
Let $\mu, \tau \in \mc{M}_1$. Then
\[
G_{\mu \boxplus (\tau + \eps \nu)}(z) = G_{\mu \boxplus \nu}(z) + \eps G_\nu(F_{\tau, \mu \boxplus \tau}(z)) F_{\tau, \mu \boxplus \tau}'(z) + o(\eps),
\]
where $F_{\tau, \mu \boxplus \tau} = K_\tau \circ G_{\mu \boxplus \tau}$. That is, $G_{D_\tau C_\mu(\nu)}(z) =  G_\nu(F_{\tau, \mu \boxplus \tau}(z)) F_{\tau, \mu \boxplus \tau}'(z)$. In particular, the operator $D_\tau C_\mu$ on $\mc{M}_0$ is well defined.
\item
Let $\set{\mu_t}$ be a free convolution semigroup. Then 
\[
G_{\mu_t \boxplus (\mu_s + \eps \nu)}(z)
= G_{s+t}(z) + \eps G_\nu(F_{s,s+t}(z)) F_{s,s+t}'(z) + o(\eps).
\]
\item
Let $\set{V_n(x,t)}$ be a family of fluctuation polynomials for $\set{\mu_t}$. Then
\[
G_{V_n^\ast(s)}(F_{s,t}(z)) F_{s,t}'(z) = G_{V_n^\ast(t)}(z).
\]
\end{enumerate}
\end{Lemma}

\begin{proof}
Clearly only the first part needs to be proven. The method of proof is similar to that of \cite[Discussion 3.5]{AnsCLT}. Denote by $\mc{I}$ the operator of functional inversion on $\mc{P}_{0,1}$. Then for $u \in \mc{P}_{0,1}$ and $v$ a formal power series with $v(0) = v'(0) = 0$, the derivative of $\mc{I}$ at $u$ in the direction $v$ is $D_u \mc{I}(v) = - (u^{-1})' v(u^{-1})$. Therefore
\begin{align*}
G_{\tau + \eps \nu}(z) &= G_\tau(z) + \eps G_\nu(z), \\
K_{\tau + \eps \nu}(z) &= K_\tau(z) - \eps K_\tau'(z) G_\nu(K_\tau(z)) + o(\eps), \\
K_{\mu \boxplus (\tau + \eps \nu)}(z) &= R_\mu + K_\tau(z) - \eps K_\tau'(z) G_\nu(K_\tau(z)) + o(\eps) \\
&= K_{\mu \boxplus \tau}(z) - \eps K_\tau'(z) G_\nu(K_\tau(z)) + o(\eps), \\
G_{\mu \boxplus (\tau + \eps \nu)}(z) &= G_{\mu \boxplus \tau}(z) + \eps G_{\mu \boxplus \tau}'(z) K_\tau'(G_{\mu \boxplus \tau}(z)) G_\nu(K_\tau(G_{\mu \boxplus \tau}(z))) + o(\eps) \\
&= G_{\mu \boxplus \tau}(z) +  \eps G_\nu(F_{\tau, \mu \boxplus \tau}(z)) F_{\tau, \mu \boxplus \tau}'(z) + o(\eps).
\end{align*}

In order to prove that the operator $D_\tau C_\mu$ on $\mc{M}_0$ is well defined, it suffices to show that for $\nu \in \mc{M}_0$, $\lim_{\eps \rightarrow 0} \frac{1}{\eps} m_n(C_\mu(\tau + \eps \nu) - C_\mu(\tau))$ exists for all $n$, since we are considering $\mc{M}_0$ with the weak$^\ast$ topology. But we have just shown that the formal generating function of these moments converges to $G_\nu(F_{\tau, \mu \boxplus \tau}(z)) F_{\tau, \mu \boxplus \tau}'(z)$.
\end{proof}

\begin{Ex}
If $\set{\mu_t}$ is a free Meixner system and $V_n(x,t) = T_n(x - \alpha(t), \beta(t))$, we can again calculate $G_{V_n^\ast(t)}$ explicitly. Indeed, $\frac{n}{2} U_n = \partial_x T_n$, so 
\[
\ip{\partial_x U_n^\ast(t)}{T_k(\cdot, t)} 
= - \ip{U_n^\ast(t)}{\partial_x T_k(\cdot, t)} 
= - \frac{n}{2} \delta_{n k},
\]  
and therefore $\partial_x U_n^\ast(t) = - \frac{n}{2} T_n^\ast(t)$. Then
\begin{align*}
G_{T_n^\ast(t)}(z) 
&= - \frac{2}{n} G_{\partial_x U_n^\ast(t)}(z)
= - \frac{2}{n} \ip{\partial_x U_n^\ast(t)}{\Res_z}
= \frac{2}{n} \ip{U_n^\ast(t)}{\partial_x \Res_z} \\
&= - \frac{2}{n} \ip{U_n^\ast(t)}{\partial_z \Res_z}
= - \frac{2}{n} \partial_z G_{U_n^\ast(t)}(z).
\end{align*}
Therefore by equation~\eqref{GChebyshev2},
\[
G_{V_n^\ast(t)}(z) 
= - \frac{2}{n} \partial_z \Bigl( u(G_t(z)) G_t^n(z) \Bigr).
\]
In particular, for the semicircular semigroup $u(z) = z$, and we recover the results of \cite{AnsCLT}.
\end{Ex}

\section{Free finite operator calculus}
\label{Sec:OperatorCalculus}

The following results are either contained in or easily deduced from \cite{RotaFiniteCalculusPaper} (reprinted in \cite{RotaFiniteCalculusBook}). All the operators involved are linear.

\begin{Prop}
\label{Prop:Rota}
Let $\set{P_n}_{n=0}^\infty$ be a sequence of polynomials such that $P_n$ is a monic polynomial of degree $n$ and $P_n(0) = 0$ for $n \geq 0$. Let $H(x,z) = \sum_{n=0}^\infty \frac{1}{n!} P_n(x) z^n$ be their exponential generating function. Let $A$ be the corresponding lowering operator on $\mf{C}[x]$, determined by $A P_n = n P_{n-1}$.  Let $W$ be the corresponding umbral operator on $\mf{C}[x]$, determined by $W x^n = P_n$. The following conditions are equivalent.
\begin{enumerate}
\item
Let $\Delta: \mf{C}[x] \rightarrow \mf{C}[x] \otimes \mf{C}[x] =
\mf{C}[x,y]$ be the usual co-multiplication $\Delta(f)(x,y) = f(x+y)$.
Then
\[
\Delta W = (W \otimes W) \Delta.
\]
\item
For $n \geq 0$,
\begin{equation*}
P_n(x + y) = \sum_{k=0}^n \binom{n}{k} P_k(x) P_{n-k}(y).
\end{equation*}
\item
$H(x + y, z) = H(x,z) H(y,z)$.
\item
$H(x,z) = e^{u(z) x}$, for some formal power series $u \in \mc{P}_{0,1}$.
\item
$\partial_x H(x,z) = u(z) H(x,z)$.
\item
$A = u^{-1}(\partial_x)$, with the same $u$.
\item
$A$ commutes with $\partial_x$.
\item
$A$ is translation invariant.
\end{enumerate}
If these conditions are satisfied, the sequence $\set{P_n}$ is called a binomial sequence.
\end{Prop}

The authors also define a Sheffer sequence to be a family of polynomials $\set{S_n(x)}$ such that $S_n(x + y) = \sum_{k=0}^n \binom{n}{k} S_k(x) P_{n-k}(y)$ for some binomial sequence $\set{P_n}$, and a cross-sequence to be a family of polynomials $\set{P_n(x,t)}$ such that $P_n(x + y, s + t) = \sum_{k=0}^n \binom{n}{k} P_k(x, s) P_{n-k}(y, t)$. They show that a large number of identities involving the classical orthogonal polynomials, and therefore a large number of combinatorial identities, are consequences of general identities for binomial sequences and cross-sequences. These results belong to the domain of umbral calculus; see \cite{dBuLoeUmbralSurvey} for a comprehensive review. The Sheffer sequences have been extensively studied by other methods as well.

It is easy to see that cross-sequences are precisely the classical Sheffer systems in the algebraic context, that is, systems of polynomials corresponding not just to classical convolution semigroups of probability measures but more generally to convolution semigroups of functionals in $\mc{M}_1$; cf.\ Remark~\ref{Rem:AllInfinitelyDivis}.

The results in the preceding proposition hold in much greater generality. All of the results in the following proposition are known, see various references in \cite{dBuLoeUmbralSurvey}. Again, all the operators involved are linear.

\begin{Prop}
\label{Prop:General}
Let $\beta$ be the sequence of positive real numbers, $\beta = \set{[0]_\beta, [1]_\beta, [2]_\beta, \ldots}$, with $[0]_\beta = 0, [1]_\beta = 1$. Denote $[n]_\beta ! = \prod_{i=1}^n [i]_\beta$, with $[0]_\beta! = 1$.

Define a linear operator $D_\beta$ on $\mf{C}[x]$ by $D_\beta (x^n) = [n]_\beta x^{n-1}$. Define the formal power series $\exp_\beta(z) = \sum_{n=0}^\infty \frac{1}{[n]_\beta !} z^n$. For $a \in \mf{R}$, define the operator $E^a_\beta = \exp_\beta(a D_\beta)$. Let $\Delta_\beta$ be the map $\mf{C}[x] \rightarrow \mf{C}[x] \otimes \mf{C}[x]$ determined by $\Delta_\beta(f)(x,a) = E^a_\beta(f)(x)$. Then $\Delta_\beta$ is co-associative.

Let $\set{P_n}_{n=0}^\infty$ be a sequence of polynomials such that $P_n$ is a monic polynomial of degree $n$ and $P_n(0) = 0$ for $n \geq 1$. Let $H(x,z) = \sum_{n=0}^\infty \frac{1}{[n]_\beta !} P_n(x) z^n$ be their generating function. Let $A$ be the corresponding lowering operator on $\mf{C}[x]$, determined by $A P_n = [n]_\beta P_{n-1}$, with the usual convention $P_{-1} = 0$. Let $W$ be the corresponding umbral operator on $\mf{C}[x]$, determined by $W x^n = P_n$. The following conditions are equivalent.
\begin{enumerate}
\item
\label{Umbral1}
$
\Delta_\beta W = (W \otimes W) \Delta_\beta.
$
\item
\label{Binomial}
For $n \geq 0$,
\begin{equation*}
\Delta_\beta(P_n)(x,y) = \sum_{k=0}^n \frac{[n]_\beta!}{[k]_\beta! [n-k]_\beta!} P_k(x) P_{n-k}(y).
\end{equation*}
\item
\label{Generating1}
$\Delta_\beta(H) = H \otimes H$.
\item
\label{ExactGenerating}
$H(x,z) = \exp_\beta(u(z) x)$, for some formal power series $u \in \mc{P}_{0,1}$.
\item
\label{DH}
$D_\beta H(x,z) = u(z) H(x,z)$.
\item
\label{Qu}
$A = u^{-1}(D_\beta)$, for the same $u$.
\item
\label{Derivative}
$A$ commutes with $D_\beta$.
\item
\label{Translation}
$A$ commutes with $E^a_\beta$ for all $a$.
\item
\label{DeltaA}
$(A \otimes I) \Delta_\beta = \Delta_\beta A$.
\end{enumerate}
\end{Prop}

\begin{proof}
$\Delta_\beta$ is co-associative, i.e. $(\Delta_\beta \otimes I) \Delta_\beta = (I \otimes \Delta_\beta) \Delta_\beta$, since both sides of this equation applied to $x^n$ give
\[
\sum_{\substack{k, l, m \geq 0 \\ k + l + m = n}} \frac{[n]_\beta!}{[k]_\beta! [l]_\beta! [m]_\beta!} x^k \otimes x^l \otimes x^m.
\]

\eqref{Umbral1} $\Leftrightarrow$ \eqref{Binomial}
$\Delta_\beta(x^n) = \sum_{k=0}^n \frac{[n]_\beta!}{[k]_\beta! [n-k]_\beta!} x^k y^{n-k}$. Thus $\Delta_\beta W x^n = \Delta_\beta P_n$, while
\[
(W \otimes W) \Delta_\beta x^n 
= (W \otimes W) \sum_{k=0}^n \frac{[n]_\beta!}{[k]_\beta! [n-k]_\beta!} x^k y^{n-k} 
= \sum_{k=0}^n \frac{[n]_\beta!}{[k]_\beta! [n-k]_\beta!} P_k(x) P_{n-k}(y).
\]

\eqref{Binomial} $\Leftrightarrow$ \eqref{Generating1}
Obvious.

\eqref{ExactGenerating} $\Leftrightarrow$ \eqref{DH}
One direction is obvious. For the other, suppose $H(x,z) = \sum_{n=0}^\infty \frac{1}{[n]_\beta!} F_n(z) x^n$, for some formal power series $\set{F_n}$. Then
\[
D_\beta H(x,z) 
= \sum_{n=1}^\infty \frac{1}{[n-1]_\beta!} F_n(z) x^{n-1} 
= \sum_{n=0}^\infty \frac{1}{[n]_\beta!} F_{n+1}(z) x^{n}.
\]
This equals $u(z) H(x,z)$ if $F_n(z) = u(z)^n$ and $H(x,z) = \exp_\beta(u(z) x)$. 

\eqref{ExactGenerating} $\Rightarrow$ \eqref{Qu}. 
The operator $A$ is determined by the equation $(A H)(x,z) = z H(x,z)$. Let $H(x,z) = \exp_\beta(u(z) x)$, with $u^{-1}(z) = \sum_{k=1}^\infty b_k z^k$. Then
\begin{align*}
u^{-1}(D_\beta) H(x,z)
&= u^{-1}(D) \sum_{n=0}^\infty \frac{1}{[n]_\beta !} u(z)^n x^n
= \sum_{n=0}^\infty \frac{1}{[n]_\beta !} u(z)^n \sum_{k=1}^\infty b_k \frac{[n]_\beta !}{[n-k]_\beta !} x^{n-k} \\
&= \sum_{m=0}^\infty \sum_{k=1}^\infty \frac{1}{[m]_\beta !}b_k u(z)^{m+k} x^m
= \sum_{m=0}^\infty u^{-1}(u(z)) \frac{1}{[m]_\beta !}u(z)^m x^m \\
&= z H(x,z),
\end{align*}
so $A = u^{-1}(D_\beta)$.

\eqref{Qu} $\Rightarrow$ \eqref{DH}
Suppose $A = u^{-1}(D_\beta)$, and so $D_\beta = u(A)$. Let $u(z) = \sum_{k=1}^\infty a_k z^k$. Then
\begin{align*}
D_\beta H(x,z)
&= u(A) H(x,z)
= u(A) \sum_{n=0}^\infty \frac{1}{[n]_\beta !} P_n(x) z^n \\
&= \sum_{n=0}^\infty z^n \sum_{k=1}^n a_k \frac{1}{[n-k]_\beta !} P_{n-k}(x)
= \sum_{m=0}^\infty \sum_{k=1}^\infty \frac{1}{[m]_\beta !} z^{m+k} a_k P_m(x) \\
&= u(z) H(x,z).
\end{align*}

\eqref{Qu} $\Rightarrow$ \eqref{Derivative}.
Obvious.

\eqref{Derivative} $\Rightarrow$ \eqref{DH}.
Suppose $A$ commutes with $D_\beta$. Then
\[
(A D_\beta H)(x,z) = D_\beta(z H(x,z)) = z (D_\beta H)(x,z).
\]
If $(D_\beta H)(x,z) = \sum_{n=1}^\infty \sum_{k=0}^{n-1} \frac{1}{[k]_\beta!} a_{n, k} P_k(x) z^n$,
\[
(A D_\beta H) (x,z) 
= \sum_{n=2}^\infty \sum_{k=1}^{n-1} \frac{1}{[k-1]_\beta!} a_{n, k} P_{k-1} (x) z^n 
= z \sum_{n=1}^\infty \sum_{k=0}^{n-1} \frac{1}{[k]_\beta!} a_{n+1, k+1} P_{k} (x) z^n.
\]
This is equal to $z (D_\beta H)(x,z)$ iff $a_{n, k} = a_{n+1, k+1}$ for $n \geq 1, k \geq 0$. But in that case, denoting $a_{n-k} = a_{n, k}$,
\[
(D_\beta H)(x,z) 
= \sum_{n=1}^\infty \sum_{k=0}^{n-1} \frac{1}{[k]_\beta!} a_{n-k} P_k(x) z^n
= \sum_{k=0}^\infty \sum_{m=1}^\infty \frac{1}{[k]_\beta!} a_m P_k(x) z^{m+k}
= u(z) H(x,z),
\]
with $u(z) = \sum_{m=1}^\infty a_m z^m$. Since $P_n$ is a monic polynomial, $a_1 = 1$.

\eqref{Derivative} $\Leftrightarrow$ \eqref{Translation}. 
This follows from the definition of $E^a_\beta$ and the relation $(D_\beta f) = \lim_{a \rightarrow 0} \frac{E^a_\beta(f) - f}{a}$.

\eqref{Translation} $\Leftrightarrow$ \eqref{DeltaA} Obvious.

\eqref{Generating1} $\Rightarrow$ \eqref{DeltaA}
$(A \otimes I) \Delta_\beta H = (A \otimes I) (H \otimes H) = z H \otimes H = z \Delta_\beta H = \Delta_\beta A H$.

\eqref{ExactGenerating} $\Rightarrow$ \eqref{Generating1}
\[
\Delta_\beta H
= \Delta_\beta \Bigl(\sum_{n=0}^\infty \frac{1}{[n]_\beta!} u(z)^n x^n \Bigl)
= \sum_{n=0}^\infty \sum_{k=0}^n \frac{1}{[k]_\beta! [n-k]_\beta!} u(z)^n x^k y^{n-k}
= H(x,z) H(y,z).
\]
\end{proof}

Most of the results from \cite{RotaFiniteCalculusPaper} have analogs in the context of Proposition~\ref{Prop:General}. We list two examples. The proofs follow from that proposition and the relation $A W = W D_\beta$ between a lowering operator and the corresponding umbral operator.

\begin{Lemma}
Call a lowering operator and an umbral operator operators appearing in these roles for some family of polynomials $\set{P_n}$ in the context of Proposition~\ref{Prop:General}.
\begin{enumerate}
\item
The umbral operators form a group, anti-isomorphic to the group $(\mc{P}_{0,1}, \circ)$. More specifically, for two umbral operators $W_1$, $W_2$ and the corresponding lowering operators $v_1(D_\beta)$, $v_2(D_\beta)$, the operator $W_2 W_1$ is the umbral operator corresponding to the lowering operator $v_1(v_2(D_\beta))$.
\item For a family of polynomials with the generating function $\exp_\beta(u(z) x)$ with $u \in \mc{P}_{0,1}$,
\[
f(u(x)) = \sum_{n=0}^\infty \frac{1}{[n]_\beta!} x^n (P_n(D_\beta)(f))(0).
\]
\end{enumerate}
\end{Lemma}

The free probability case corresponds to the sequence $\beta$ with $[n]_\beta = 1$ for $n \geq 1$. This case has already been singled out, for example, in \cite[Part XII]{JoniRota}. Namely, the classical case $[n]_\beta = n$ is distinguished by the fact that $\Delta_\beta$ is a homomorphism, and so $\mf{C}[x]$ with the usual multiplication and co-multiplication $\Delta_\beta$ is a bialgebra. The free probability case $[n]_\beta = 1$ is distinguished by the fact that $\partial_\beta$, where $\partial_\beta(f)(x,y) = \frac{\Delta_\beta(f)(x,y) - f(y)}{x}$, is a derivation, and so $\mf{C}[x]$ with the usual multiplication and co-multiplication $\partial_\beta$ is an example of an \emph{infinitesimal coalgebra} of \cite{JoniRota} or of a \emph{generalized difference quotient ring} of \cite{VoiCoalgebra}. Note that the classical case is also the only one when $D_\beta$ is a derivation.

In the free case, in addition to the Proposition \ref{Prop:General} being valid, we also have \emph{different} analogs of the first three statements of Proposition \ref{Prop:Rota}. Namely, in this case $\mf{C}[x]$, in addition to being a coalgebra with co-multiplication $\Delta_\beta$, also has the structure of a dual group in the sense of \cite{VoiDualGroup} with an operation $\Delta_\ast$. In the following proposition, for a pair of algebras $\mc{A}_1, \mc{A}_2$, by $\mc{A}_1 \ast \mc{A}_2$  we mean the algebraic reduced free product of algebras (amalgamated over the identity element). Clearly the free product $\mf{C}[x] \ast \mf{C}[x]$ is isomorphic to $\mf{C}\langle x,y\rangle$, the algebra of polynomials in two non-commuting variables. For two unital operators $W_1$ on $\mc{A}_1$, $W_2$ on $\mc{A}_2$, their free product operator is the unital operator on $\mc{A}_1 \ast \mc{A}_2$ determined by $(W_1 \ast W_2) (a_1 a_2 \cdots a_n) = (W_{i(1)} a_1) (W_{i(2)} a_2) \cdots (W_{i(n)} a_n)$, where $a_j \in \mc{A}_{i(j)} \backslash \mf{C}$ and for all $j$, $i(j) \neq i(j+1)$.

\begin{Prop}
\label{Prop:FreeBinomial}
Let $\beta$ be a sequence with $[n]_\beta = 1$ for $n \geq 1$. Consider the corresponding objects from Proposition \ref{Prop:General}. More specifically, in this case
\begin{gather*}
D_\beta f = \frac{f(x) - f(0)}{x} = (\partial f)(x,0), \qquad
\exp_\beta(z) = \frac{1}{1-z}, \\
E^a_\beta(f)(x) = \frac{x f(x) - a f(a)}{x-a} = x (\partial f)(x,a) + f(a), \\
\Delta_\beta(f)(x,y) = \frac{x f(x) - y f(y)}{x-y} = x (\partial f)(x,y) + f(y), \\
H(x,z) = \sum_{n=0}^\infty P_n(x) z^n, \qquad
A P_n = P_{n-1}.
\end{gather*}
Then the following conditions are equivalent to the conditions of Proposition \ref{Prop:General} for this $\beta$.
\begin{enumerate}
\item
\label{Umbral}
Let $\Delta_\ast: \mf{C}[x] \rightarrow \mf{C}[x] \ast \mf{C}[x] =
\mf{C}\langle x,y\rangle$ be the homomorphism $\Delta_\ast(f)(x,y) = f(x+y)$.
Then
\begin{equation}
\label{UmbralEquation}
\Delta_\ast W = (W \ast W) \Delta_\ast.
\end{equation}
\item
\label{Free}
For $n \geq 1$,
\begin{multline}
\label{FreeEquation}
P_n(x+y) = \sum_{k=1}^n \sum_{\substack{i(1), i(2), \ldots, i(k) \geq 1 \\ i(1) + i(2) + \ldots + i(k) = n}} (P_{i(1)}(x) P_{i(2)}(y) P_{i(3)}(x) P_{i(4)}(y) \cdots \\
+ P_{i(1)}(y) P_{i(2)}(x) P_{i(3)}(y) P_{i(4)}(x) \cdots).
\end{multline}
\item
\label{Generating}
The generating function $H$ satisfies the equation
\[
H(x + y,z) = \frac{H(x,z) H(y,z)}{H(x,z) + H(y,z) - H(x,z) H(y,z)}.
\]
\item
\label{ExactGenerating1}
$H(x,z) = \frac{1}{1 - u(z) x}$ for some formal power series $u \in \mc{P}_{0,1}$.
\end{enumerate}
\end{Prop}

\begin{proof}
\eqref{Umbral} $\Leftrightarrow$ \eqref{Free}. 
Obvious.

\eqref{Free} $\Leftrightarrow$ \eqref{Generating}. 
Equation ~\eqref{FreeEquation} is equivalent to
\begin{align*}
H(x+y, z) &= 1 + \sum_{n=1}^\infty \sum_{k=1}^n \sum_{\substack{i(1), i(2), \ldots, i(k) \geq 1 \\ i(1) + i(2) + \ldots + i(k) = n}} \Bigl( P_{i(1)}(x) z^{i(1)} P_{i(2)}(y) z^{i(2)} P_{i(3)}(x) z^{i(3)} P_{i(4)}(y) z^{i(4)} \cdots \\*
&\quad + P_{i(1)}(y) z^{i(1)} P_{i(2)}(x) z^{i(2)} P_{i(3)}(y) z^{i(3)} P_{i(4)}(x) z^{i(4)} \cdots \Bigr) \\
&= 1 + \sum_{k=1}^\infty \biggl( \underbrace{(H(x,z) - 1) (H(y,z) - 1) (H(x,z) - 1) \cdots}_{k \text{ terms}} \\*
&\quad + \underbrace{(H(y,z) - 1) (H(x,z) - 1) (H(y,z) - 1) \cdots}_{k \text{ terms}} \biggr) \\
&= H(x,z) \frac{1}{1 - (H(y,z) - 1) (H(x,z) - 1)} H(y,z) \\
&= \frac{H(x,z) H(y,z)}{H(x,z) + H(y,z) - H(x,z) H(y,z)}.
\end{align*}

\eqref{Generating} $\Leftrightarrow$ \eqref{ExactGenerating1}.
 This follows from the fact that the equation
\[
\frac{1}{H(x+y)} = \frac{1}{H(x)} + \frac{1}{H(y)} - 1
\]
has a unique formal power series solution $\frac{1}{H(x)} = 1 - a x$, for $a \in \mf{C}$. Therefore $H(x,z) = \frac{1}{1 - u(z) x}$. Since $H(x,0) = 1$ and $P_n$ is a monic polynomial of degree exactly $n$, we conclude that $u(0) = 0$, $u'(0) = 1$.
\end{proof}

Note that the generating functions of the (algebraic) free Sheffer systems are precisely those satisfying 
\[
H(x + y, s + t, z) = \frac{H(x,s,z) H(y,t,z)}{H(x,s,z) + H(y,t,z) - H(x,s,z) H(y,t,z)}.
\]
Indeed, by the same argument as before, the solutions of this equation are precisely of the form $\frac{1}{H(x,t)} = 1 + a x + b t$, for $a, b \in \mf{C}$. Therefore $H(x,t,z) = \frac{1}{1 + t f(z) - u(z) x}$. Since $P_n$ is a monic polynomial of degree exactly $n$, $u(0) = 0$, $u'(0) = 1$, $f(0) = 0$. Therefore we can define $R_\mu$ by $R_\mu(z) = \frac{1}{z} f(u^{-1}(z))$, so that
\[
H(x,t,z) = \frac{1}{u(z) (K_t(u(z)) - x)} = \frac{1}{1 - u(z)(x - t R_\mu(u(z)))}.
\]
Note that $\mu \in \mc{M}_1$, it need not be a positive measure. We leave to the reader the analogs of equations~\eqref{UmbralEquation} and \eqref{FreeEquation} for the free Sheffer systems. Cf.\ \cite{LehSpectra}, where such expansions were considered in the more general operator-valued case.

\begin{Ex}[Free Meixner systems]
Time-zero polynomials of a free Sheffer system satisfy the conditions of Proposition~\ref{Prop:FreeBinomial}. For the free Meixner systems, these are the following.
\begin{description}
\item[Hermite (Chebyshev)]
$P_n(x) = x^n$. 
\item[Charlier]
$P_0(x) = 1$, $P_n(x) = x (x-1)^{n-1}$ for $n \geq 1$. 
\item[Meixner / Laguerre / Meixner-Pollaczek] $P_0(x) = 1$, $P_n(x) = x U_{n-1}(x - 2a)$ for $n \geq 1$. Indeed, for the modified Chebyshev polynomials, the time-zero polynomials are $Q_n(x) = U_n(x - 2a)$. Thus from equations~\eqref{PQ}, the time-zero polynomials of the free Meixner system with parameter $a$ are $P_n(x) = U_n(x - 2a) + 2a U_{n-1}(x - 2a) + U_{n-2}(x - 2a) = x U_{n-1}(x - 2a)$.
\end{description}
\end{Ex}

The sequence $P_n(x) = x^n$ is a binomial sequence for any $\beta$. However, already for the zero-time polynomials of the ``$\beta$-Charlier systems'', this is in general not the case.

\begin{Lemma}
\label{Lemma:Counter}
The sequence $P_n(x) = \prod_{k=0}^{n-1} (x - [k]_\beta)$ is a $\beta$-binomial sequence in only two cases: $[n]_\beta = n$, or $[n]_\beta = 1$ for $n \geq 1$.
\end{Lemma}

\begin{proof}
If the sequence $P_n(x) = \prod_{k=0}^{n-1} (x - [k]_\beta)$ is a $\beta$-binomial sequence, then the corresponding lowering operator, determined by $A P_n = [n]_\beta P_{n-1}$, commutes with $D_\beta$. Evaluating the identity $A D_\beta = D_\beta A$ on $x^n$ and comparing coefficients, we obtain a system of equations in $\set{[j]_\beta}$ that is linear in $[n]_\beta$. So by induction, the sequence $\set{[n]_\beta}$ is determined by $[2]_\beta$. From this analysis for $x^3$, we conclude that $[3]_\beta = \frac{1}{2} ([2]_\beta + [2]_\beta^2)$, and from this analysis for $x^4$ we conclude that $[2]_\beta$ is a root of the equation $(x-1)^2 (x-2) (x+1)$. Since $[2]_\beta \geq 0$, it is equal to either $2$ or $1$. In the first case, $[n]_\beta = n$, which corresponds to the classical situation. In the second case, $[n]_\beta = 1$ for $n \geq 1$, which corresponds to the free situation.
\end{proof}

Moreover, already for $q$-Hermite polynomials, their generating function is not of the form appearing in Proposition \ref{Prop:General}. In particular, it is not clear if the families of Remark~\ref{Rem:QDeform} can be interpreted as Sheffer systems, or even if their time-zero polynomials can be interpreted as binomial systems. On the other hand, some positive evidence for the Sheffer interpretation is provided by the fact \cite{BKSQGauss} that the $q$-Hermite polynomials are martingale polynomials for the $q$-Brownian motion with respect to its standard filtration. Also, we have recently learned from Professor Ismail about a different deformation of the umbral calculus, which may be appropriate for these families. This direction will be pursued in a future paper.


\providecommand{\bysame}{\leavevmode\hbox to3em{\hrulefill}\thinspace}

\end{document}